# Numerical-Analytical Investigation into Impact Pipe Driving in Soil with Dry Friction. Part I: Nondeformable External Medium

**N. I. Aleksandrova**

*N.A. Chinakal Institute of Mining, Siberian Branch, Russian Academy of Sciences, Krasnyi pr. 54, Novosibirsk, 630091 Russia*
*e-mail: alex@math.nsc.ru*

**Abstract**—The study focuses on propagation of longitudinal waves in an elastic pipe partly embedded in a medium with dry friction. Mathematical formulation of the problem on the impact pipe driving into the soil is based on the model of longitudinal vibration of an elastic rod with taking into account lateral resistance. The lateral resistance of soil is described by the law of the contact dry friction. Numerical and analytical solutions to problems on longitudinal impulse loading of a pipe are compared.

## 1. INTRODUCTION

The problems of pipe "behavior" in the soil arise in a variety of technological processes including the impact driving and removal of piles, the trenchless underground communication pipe laying by driving metallic pipes into the soil, behavior of underground pipelines in earthquakes, or different rod elements in mechanical systems. The influence of the friction between the soil and the side pipe surface on the wave process anyway remains a research issue of prime importance. A good deal of R&D reports and transactions reviews the interaction of bodies with account for the friction factor [1–21].

The problem of the longitudinal wave propagation in a rod which surface interacts with an ambient medium like a dry-friction model, was studied in [3, 4, 6, 7, 11–21] by considering and analyzing the wave processes. In most publications the researchers established a rod cross-section, partitioning a disturbed section from a stationary section of the rod. The general method for solving these problems implies the reduction of the study non-linear problem to a number of soluble consistently linear problems. The basic problem the researchers face is the identification of the range boundaries in time and in the longitudinal coordinate within which a magnitude and direction of the friction force remain constant. A friction force direction depends on a velocity sign in these waves. The solution in every range is performed with consideration for a preliminary determined sign of the velocity. The approach used in [3, 4, 6, 7, 11–18] made it possible to obtain a solution to a number of the problems under rather simple rod-loading laws.

In [3] the method of characteristics was helpful to study the problem of an impact effect against an elastic rod, partially embedded into a dry-friction medium. The problems of an impact at a constant rate and a rate monotonously increasing in time are solved. The researchers found the critical impact rate at which the impact impulse front reaches the opposite end of the rod.

In [4] the study subject is a moving rod of a finite length with account for a contact dry friction between a rod and an ambient medium. A constant motion velocity is preset at the left end of the rod. A motion of the ambient medium is taken into account solely in the velocity of the moving lateral load applied. In this problem statement the analytical solutions are obtained by using Laplace transformations.

In [7] the research problem deals with the disturbance propagation within a rod at the side surface of which the dry-friction law is realized when one of rod ends moves according to the harmonic law. The approximated solutions were gained by the methods of the integrated amplitudes and the mechanical analogue. These two solutions coincide qualitatively and still differ quantitatively.

In [6, 11, 15] the researchers obtained an analytical solution to the problem of exciting oscillations in a rod with the external dry friction and one rod end being under the constant stress and the opposite rod end stiffly fixed. In [12 – 16] under consideration is the problem dealing with the dynamics of an elastic rod, one of its ends interacts with a mass striking it and at the other end either a displacement or stress-displacement relation is preset pertaining to a character of the load applications. Besides the gravity force of the rod and the stiff body, as well as the resistance under the law of the dry friction at the rod side surface were taken into account. In [16] the solution to this problem is reported for a random time moment and a precised time of a weight rebound from the rod.

In [17, 18] the studies cover the longitudinal oscillations of the finite elastic rod with a preset dry friction at the side surface and the stepwise loading of one of its ends with a period, being a multiple of a natural oscillation period for the system. The resonance excitation and the stabilized periodical oscillations of the system are analyzed.

Despite the great theoretical contribution gained in solving the problem the issues of the prime practical importance are left unsolved. For example, there are significant constraints on a character of the load applied to the rod, because under the smooth loading the front of the wave propagating along the rod moves slower as compared to a disturbance propagation velocity in the rod and it can be evaluated exclusively in solutions, but the latter is not always possible. Thereto, the most practical problems of the loading are complex.

Publications [19, 21] present the numerical solution. In [19] the solution to the problem of the wave propagation in a viscous –elastic rod with the external dry friction under a rectangular dynamic effect was obtained by the method of characteristics. In [21] the problem of the dynamic effect on a pipe with the Coulomb’s dry friction was solved by the finite difference method.

Below the similar to [21] method for the end-to-end numerical calculation of the wave processes in the systems with dry friction permits to consider any load type because a direction and a magnitude of the friction force at every point and every time moment is preset in terms of the physical aspects in the course of a solving procedure. In this case the main problem is the testing of numerical algorithms and achievement of adequacy between discrete and continuum solutions. In view of this point the availability of analytical solutions is a compulsory condition. On the other hand, the analytical solutions in terms of some simplifying assumptions should be also verified on the correctness and the identification of the parameter domain of their validity.

In the present paper the research scientist obtained the approximated analytical solutions to the problem of the disturbance propagation in a rod of a finite and infinite lengths with the external dry friction when a pulse stress of a random form is preset at one of its ends, while the other end is free of stresses.

## 2. STATEMENT OF THE PROBLEM

In the study problem an elastic tubular rod with radius $R$, wall thickness $h$, and length $L$ is embedded in the soil at depth $L_1$. A longitudinal pulse $Q(t)$ is applied to its end. The constant shear stress $\tau_0$ is applied in moving cross-sections at a side surface of the rod. The dry friction effect is accounted in the way similar to [13]. The origin of the selected coordinate system should coincide with the rod end, subjected to a percussion effect, axis $Z$ should be directed in parallel to the rod axis in-depth the medium. The rod motion is described by 1D wave equation relative to displacements $U(t,z)$:

$$\ddot{U} = c^2 U''_{,zz} - k\tau\,, \quad \tau = \frac{\tau_0 P_t}{\rho S_t}\,, \quad k = \text{sign}\dot{U}\,, \tag{1}$$

under zero initial conditions:

$$U\big|_{t=0} = 0\,, \quad \dot{U}\big|_{t=0} = 0 \tag{2}$$

and boundary conditions at the rod ends:

$$ES_t U'_{,z}\big|_{z=0} = -Q(t)\,, \quad ES_t U'_{,z}\big|_{z=L} = 0\,. \tag{3}$$

In (1) –(3) $c = \sqrt{E/\rho}$ is a velocity of a longitudinal wave, $\rho$ is density of a pipe material, $E$ is Young modulus, $S_t = \pi h(2R - h)$ is area of the pipe cross-section, $P_t = 2\pi R$ is pipe perimeter.

## 3. SPECIFIC FEATURES OF NUMERICAL ALGORITHMS

The system of equations (1)–(3) was solved numerically for a single percussion. The finite difference scheme of “cross” type was used:

$$U_j^{n+1} - 2U_j^n + U_j^{n-1} = h_t^2\left[c^2(U_{j+1}^n - 2U_j^n + U_{j-1}^n)\big/h_z^2 - k\tau_j^n\right].$$

Here $U_j^n = U(h_t n, h_z j)$ are displacements at time moments $t = h_t n$ at point coordinates $z = h_z j$, where $h_t$, $h_z$ are mesh of the difference net in coordinates $t$, $z$, $\tau_j^n = \tau(h_t n, h_z j)$. The optimal parameters of the difference net to provide the minimum numerical dispersion for a discrete analogue of equation (1) are $h_z = ch_t$, where $h_t$, $h_z$ are mesh of the difference net in coordinates $t$, $z$ [22].

The calculation algorithm with account for the friction is: as neither direction nor a friction force value are known the velocities of points for two possible signs are calculated at $k$ ($k > 0$ and $k < 0$):

a) in the first case introduce a fictive velocity $\dot{U}_j^+$:

$$\dot{U}_j^+ = \frac{U_j^+ - U_j^n}{h_t}\,, \quad U_j^+ = U_j^{n+1} - h_t^2\tau\,;$$

b) in the second case introduce a fictive velocity $\dot{U}_j^-$:

$$\dot{U}_j^- = \frac{U_j^- - U_j^n}{h_t}\,, \quad U_j^- = U_j^{n+1} + h_t^2\tau\,.$$

In the above two cases $U_j^{n+1}$ is calculated from the finite difference analogue of equation (1) disregarding the friction.

Two situations can be observed:

1) at $\dot{U}_j^+$ and $\dot{U}_j^-$ of the same sign the value $U_j^k$ for which $\text{abs}(\dot{U}_j^k) = \min[\text{abs}(\dot{U}_j^-), \text{abs}(\dot{U}_j^+)]$ is selected as the real values of displacements $U_j^{n+1}$ from $U_j^+$ and $U_j^-$;

2) at $\dot{U}_j^+$ and $\dot{U}_j^-$ of different signs or one of the velocities turned to zero and under assumption of the passive friction it is concluded that the real velocity of the pipe equals to zero and, therefore, the friction force is absent. Values $U_j^{n+1}$ are calculated at assumed $k=0$.

Thus, the identification of boundaries, separating the motion and rest areas, being one of main hampering problems in the analytical solution, is reduced to finding points where $\dot{U}_j^+$ and $\dot{U}_j^-$ of different signs or one of them is turned to zero. As a magnitude and a direction of the friction force are single-valued in calculations, every temporal layer implies solving a linear problem where the friction force is found and expressed as a load in the right side of the equation.

## 4. ANALYTICAL ESTIMATES

To test the finite-difference algorithms we obtain the approximated analytical estimates for a number of problems. Suggest that a semi-sinusoidal pulse load is applied to a section $z=0$ of a semi-infinite rod with the acting dry contact friction at its side surface:

$$Q(t)=P_0\sin(\omega_* t)H_0(t_0-t)H_0(t),\quad \omega_*=\pi/t_0. \tag{4}$$

To obtain an analytical solution we assume that the rod is infinite and a load is applied to cross-section $z=0$ and perform Laplace transformation in time and Fourier transformation along the longitudinal coordinate $z$. As the velocity in the infinite rod is positive, then $k=1$ if $\dot{U}>0$, and $k=0$ if $\dot{U}=0$. Note that the disturbance propagation velocity in a waveguide with account for the friction can not exceed a velocity value $c$. In view of these two remarks the solution in the images of the velocity of the pipe section displacements takes the form:

$$\dot{U}^{LF}=\frac{P_0(1+e^{-pt_0})\omega_* p}{\rho S_t(\omega_*^2+p^2)(p^2+q^2c^2)}-\frac{2\tau_0 P_t cpk}{\rho S_t(p^2+q^2c^2)^2}.$$

Conversion of Laplace and Fourier transformations gives the following displacement velocity expression for an infinite rod:

$$\dot{U}(z,t)=\frac{1}{c\rho S_t}\left[\frac{P_0}{2}\sin\omega_*\left(t-\frac{|z|}{c}\right)H_0\left(\frac{|z|}{c}-t+t_0\right)-k\frac{\tau_0 P_t ct}{2}\right]H_0\left(t-\frac{|z|}{c}\right).$$

In the case with a semi-infinite rod $(z\geq 0)$ we should provide a twofold increase in the pulse amplitude:

$$\dot{U}(z,t)=\frac{1}{c\rho S_t}\left[P_0\sin\omega_*\left(t-\frac{z}{c}\right)H_0\left(\frac{z}{c}-t+t_0\right)-k\frac{\tau_0 P_t ct}{2}\right]H_0\left(t-\frac{z}{c}\right).$$

Given that $\dot{U}>0$, the function $k=1$, so the result is an inequality for identification of a range at the plane $(z,\, t)$ where $k=1$:

$$\sin\omega_*\left(t-\frac{z}{c}\right)>\frac{\tau_0 P_t ct}{2P_0},\quad 0\leq t-\frac{z}{c}\leq t_0.$$

For the given time moment $t$ we find the interval in $z$, where this inequality is fulfilled. Let $z_1$ satisfies the equation:

$$\sin\omega_*\left(t-\frac{z_1}{c}\right)=\frac{\tau_0 P_t ct}{2P_0}.$$

Considering symbols $\varepsilon = t - z_1 / c$, we obtain that $k = 1$ within the interval $c(t - t_0 + \varepsilon) < z < c(t - \varepsilon)$. Beyond this interval $\dot{U} = 0$ and $k = 0$. Therefore, we have the following approximated relation of the displacement velocity vs. coordinate $z$ at fixed time $t$:

$$\dot{U}(z,t) = \frac{1}{c\rho S_t}\left[P_0 \sin\omega_*\left(t - \frac{z}{c}\right) - \frac{\tau_0 P_t ct}{2}\right] H_0\left(\frac{z}{c} - t + t_0 - \varepsilon\right) H_0\left(t - \frac{z}{c} - \varepsilon\right). \quad (5)$$

Here $\varepsilon$ is found from the term:

$$\sin(\omega_* \varepsilon) = \frac{\tau_0 P_t ct}{2P_0} \quad (6)$$

and satisfies the inequality $0 \le \varepsilon \le t_0 / 2$.

It is revealed from the analysis of formula (5) that the maximum displacement velocity values are gained in the cross-section $z = c(t - t_0 / 2)$ and linearly depend on a pulse amplitude $P_0$ of shear stress $\tau_0$ and time $t$:

$$\max \dot{U} = \frac{P_0}{c\rho S_t} - \frac{\tau_0 P_t t}{2\rho S_t} > 0.$$

It follows from this formula that in a semi-infinite rod with an applied external dry friction a running pulse amplitude disappears at time moment:

$$t_* = \frac{2P_0}{\tau_0 P_t c}. \quad (7)$$

The estimate (7) coincides with the estimate, obtained in [6].

As the disturbance propagation velocity in the rod equals to $c$, it follows from (7) that at

$$z > z_* = c\left(t_* - \frac{t_0}{2}\right) = \frac{2P_0}{\tau_0 P_t} - \frac{ct_0}{2} \quad (8)$$

no disturbances are present in the rod.

Let fix coordinate $z$ and establish the dependence of the displacement velocity on time. Let $t_1$, $t_2$ satisfy equations:

$$\sin\omega_*\left(t_1 - \frac{z}{c}\right) = \frac{\tau_0 P_t ct_1}{2P_0}, \quad \sin\omega_*\left(t_2 - \frac{z}{c}\right) = \frac{\tau_0 P_t ct_2}{2P_0}.$$

With account for symbols $\varepsilon_1 = t_1 - z / c$, $\varepsilon_2 = z / c + t_0 - t_2$ we obtain transcendental equations for evaluation of $\varepsilon_1$, $\varepsilon_2$:

$$\sin(\omega_* \varepsilon_1) = \frac{\tau_0 P_t c}{2P_0}\left(\varepsilon_1 + \frac{z}{c}\right), \quad \sin(\omega_* \varepsilon_2) = \frac{\tau_0 P_t c}{2P_0}\left(\frac{z}{c} + t_0 - \varepsilon_2\right),$$

where $\varepsilon_1$, $\varepsilon_2$ lie in he interval $0 \le \varepsilon_1 \le t_0 / 2$, $0 \le \varepsilon_2 \le t_0 / 2$. We obtain that $\dot{U} > 0$ and $k = 1$ in the domain $z / c + \varepsilon_1 < t < z / c + t_0 - \varepsilon_2$. Beyond this interval $\dot{U} = 0$. Now let take the following relation of the displacement velocity vs. time at a fixed coordinate $z$ as an approximated solution:

$$\dot{U}(z,t) = \frac{1}{c\rho S_t}\left[P_0 \sin\omega_*\left(t - \frac{z}{c}\right) - \frac{\tau_0 P_t ct}{2}\right] H_0\left(\frac{z}{c} - t + t_0 - \varepsilon_2\right) H_0\left(t - \frac{z}{c} - \varepsilon_1\right). \quad (9)$$

Integrating (9), we obtain the dependence of the rod cross-section displacements on the time and the coordinate of the semi-sinusoidal pulse:

$$U(z,t)=\frac{1}{c\rho S_t}\begin{cases}\frac{P_0}{\omega_*}[\cos(\omega_*\varepsilon_1)+\cos(\omega_*\varepsilon_2)]-\frac{\tau_0 P_t c}{4}\left[\left(\frac{z}{c}+t_0-\varepsilon_2\right)^2-\left(\frac{z}{c}+\varepsilon_1\right)^2\right], & t-\frac{z}{c}>t_0-\varepsilon_2,\\ \frac{P_0}{\omega_*}\left[\cos(\omega_*\varepsilon_1)-\cos\omega_*\left(t-\frac{z}{c}\right)\right]-\frac{\tau_0 P_t c}{4}\left[t^2-\left(\frac{z}{c}+\varepsilon_1\right)^2\right], & \varepsilon_1\le t-\frac{z}{c}\le t_0-\varepsilon_2,\\ 0, & t-\frac{z}{c}<\varepsilon_1.\end{cases} \quad (10)$$

Given that a random pulse effect is preset in the interval $(0,t_0)$

$$Q(t)=P_0\tilde{Q}(t)H_0(t_0-\ t)H_0(\ t),$$

where $\tilde{Q}(t)$ is a convex function, we follow the analogous procedure and obtain an approximated dependence of the displacement velocity on the time $t$ at a fixed coordinate $z$:

$$\dot{U}(z,t)=\frac{1}{c\rho S_t}\left[P_0\tilde{Q}\left(t-\frac{z}{c}\right)-\frac{\tau_0 P_t ct}{2}\right]H_0\left(\frac{z}{c}-t+t_0-\varepsilon_2\right)H_0\left(t-\frac{z}{c}-\varepsilon_1\right), \quad (11)$$

$$2P_0\tilde{Q}(\varepsilon_1)=\tau_0 P_t(z+c\varepsilon_1),\quad 2P_0\tilde{Q}(t_0-\varepsilon_2)=\tau_0 P_t(z+ct_0-c\varepsilon_2).$$

In the case of a rectangular-shaped pulse $\varepsilon_1=\varepsilon_2=0$ and the displacements resulted from the integration (11), no doubt, satisfy equation (1). For other-shaped pulses equation (1) in solution (11) is satisfied approximately, but with a higher precision in the case when the term $t_0 << t_*$ is fulfilled.

Note that (5), (10), (11) are derived under assumption that $k=1$ if $\dot{U}>0$, and $k=0$ if $\dot{U}=0$. This term is fulfilled within the entire interval of the disturbance propagation under action of a correct-shaped pulse, given that $t_0\le t_*$. If $t_0>t_*$, a displacement velocity can become a negative value at a certain time moment under assumption $k=-1$. In this case the solution will be different.

Consider a rod of a finite length $L$ with a dry friction, a shear stress $\tau_0$, and semi-sinusoidal longitudinal load all applied to its side surface (4).

It follows from (9), in a rod of a length $L>z_*$ the end being opposite to a point of the force application does not "feel" this effect. It is possible to evaluate from (9) what pulse amplitude is required to enable a disturbance to reach a pipe end:

$$P_0\ge\left(L+\frac{ct_0}{2}\right)\frac{\tau_0 P_t}{2}.$$

If $L>>ct_0/2$, the second term can be neglected and the result is:

$$P_0\ge F_{\text{тр}}/2,\quad F_{\text{тр}}=\tau_0 P_t L, \quad (12)$$

viz., a pulse amplitude should at least a half higher than a friction force $F_{\text{f}}$, acting at the side surface of the pipe.

Taking into account the rebounds from the rod ends, the approximated solution to the study problem for displacement velocities is:

$$\dot{U}(z,t) = \sum_{n=0}^{n_1} \frac{1}{c\rho S_t}\left\{P_0 \sin\omega_*(t - a_{2n} + \varepsilon_{11n}) - F_{\text{тр}}\frac{ct}{2L}\right\} H_0(t-a_{2n})H_0(a_{1n}-t) + \\ + \sum_{n=0}^{n_2} \frac{1}{c\rho S_t}\left\{P_0 \sin\omega_*(t - b_{2n} + \varepsilon_{21n}) - F_{\text{тр}}\frac{ct}{2L}\right\} H_0(t-b_{2n})H_0(b_{1n}-t), \tag{13}$$

$$n_1 = \left[\frac{P_0}{F_{\text{тр}}} - \frac{z}{2L}\right], \quad n_2 = \left[\frac{P_0}{F_{\text{тр}}} + \frac{z}{2L} - 1\right],$$

$$\sin(\omega_*\varepsilon_{11n}) = \frac{F_{\text{тр}}}{P_0}\left(n + \frac{z + c\varepsilon_{11n}}{2L}\right), \quad \sin(\omega_*\varepsilon_{12n}) = \frac{F_{\text{тр}}}{P_0}\left(n + \frac{z + c(t_0 - \varepsilon_{12n})}{2L}\right),$$

$$\sin(\omega_*\varepsilon_{21n}) = \frac{F_{\text{тр}}}{P_0}\left(n + 1 - \frac{z - c\varepsilon_{21n}}{2L}\right), \quad \sin(\omega_*\varepsilon_{22n}) = \frac{F_{\text{тр}}}{P_0}\left(n + 1 - \frac{z - c(t_0 - \varepsilon_{22n})}{2L}\right),$$

$$a_{1n} = \frac{2Ln + z}{c} + t_0 - \varepsilon_{12n}, \quad a_{2n} = \frac{2Ln + z}{c} + \varepsilon_{11n},$$

$$b_{1n} = \frac{2L(n+1) - z}{c} + t_0 - \varepsilon_{22n}, \quad b_{2n} = \frac{2L(n+1) - z}{c} + \varepsilon_{21n}.$$

The square brackets in (13) denote an integral part of a number. Integration of (13) in time gives the relation of the rod cross-section displacements vs. time and a longitudinal coordinate:

$$U(z,t) = \sum_{n=0}^{n_1} U_{1n} + \sum_{n=0}^{n_2} U_{2n}, \tag{14}$$

$$U_{1n}(z,t) = \frac{1}{c\rho S_t}\begin{cases} \dfrac{P_0}{\omega_*}\left[\cos(\omega_*\varepsilon_{11n}) + \cos(\omega_*\varepsilon_{12n})\right] - \dfrac{\tau_0 P_t c}{4}(a_{1n}^2 - a_{2n}^2), & t > a_{1n}, \\ \dfrac{P_0\left[\cos(\omega_*\varepsilon_{12n}) - \cos\omega_*(t - a_{2n} + \varepsilon_{11n})\right]}{\omega_*} - \dfrac{\tau_0 P_t c}{4}(t^2 - a_{2n}^2), & a_{2n} \le t \le a_{1n}, \\ 0, & t < a_{2n}, \end{cases}$$

$$U_{2n}(z,t) = \frac{1}{c\rho S_t}\begin{cases} \dfrac{P_0}{\omega_*}\left[\cos(\omega_*\varepsilon_{21n}) + \cos(\omega_*\varepsilon_{22n})\right] - \dfrac{\tau_0 P_t c}{4}(b_{1n}^2 - b_{2n}^2), & t > b_{1n}, \\ \dfrac{P_0}{\omega_*}\left[\cos(\omega_*\varepsilon_{22n}) - \cos\omega_*(t - b_{2n} + \varepsilon_{21n})\right] - \dfrac{\tau_0 P_t c}{4}(t^2 - b_{2n}^2), & b_{2n} \le t \le b_{1n}, \\ 0, & t < b_{2n}. \end{cases}$$

Evaluate a cumulative slip of a pipe in the cross-section $z = 0$ at a time moment $t > b_{1n_*}$, namely after the disturbances in a pipe are decayed under action of the friction force. It follows from (14):

$$U_1^* = \frac{1}{c\rho S_t}\sum_{n=0}^{n_*}\left(\frac{P_0}{\omega_*}[\cos(\omega_*\varepsilon_{11n})+\cos(\omega_*\varepsilon_{12n})]-\frac{F_{\text{тр}}c}{4L}(t_0-\varepsilon_{11n}-\varepsilon_{12n})\left(t_0+\frac{4Ln}{c}+\varepsilon_{11n}-\varepsilon_{12n}\right)\right)+$$
$$+\frac{1}{c\rho S_t}\sum_{n=0}^{n_*-1}\left(\frac{P_0}{\omega_*}[\cos(\omega_*\varepsilon_{21n})+\cos(\omega_*\varepsilon_{22n})]-\frac{F_{\text{тр}}c}{4L}(t_0-\varepsilon_{21n}-\varepsilon_{22n})\left(t_0+\frac{4L(n+1)}{c}+\varepsilon_{21n}-\varepsilon_{22n}\right)\right), \quad (15)$$

$$n_* = \left[\frac{P_0}{F_{\text{тр}}}\right], \quad \sin(\omega_*\varepsilon_{11n}) = \frac{F_{\text{тр}}}{P_0}\left(n+\frac{c\varepsilon_{11n}}{2L}\right), \quad \sin(\omega_*\varepsilon_{12n}) = \frac{F_{\text{тр}}}{P_0}\left(n+\frac{c(t_0-\varepsilon_{12n})}{2L}\right),$$

$$\sin(\omega_*\varepsilon_{21n}) = \frac{F_{\text{тр}}}{P_0}\left(n+1+\frac{c\varepsilon_{21n}}{2L}\right), \quad \sin(\omega_*\varepsilon_{22n}) = \frac{F_{\text{тр}}}{P_0}\left(n+1+\frac{c(t_0-\varepsilon_{22n})}{2L}\right).$$

Assume that $n_*$ is high enough, and approximated sum can be found by integration instead of the summing. We obtain a rough estimate of the cumulative slip of the pipe as a result of multiple reflections of a semi-sinusoidal pulse:

$$U_1^* = \frac{t_0}{2c\rho S_t F_{\text{тр}}}\left\{P_0^2+2P_0F_{\text{тр}}\left(1+\frac{2}{\pi}\left(1+\frac{t_0c}{2L}\right)\right)-F_{\text{тр}}^2\frac{t_0c}{2L}\right\}. \quad (16)$$

From the formula if the impact energy equals to $t_0P_0^2/2$, the cumulative slip is proportional to "impact energy with account for friction" $t_0\tilde{P}_0^2/2$, where

$$\tilde{P}_0^2 = \left\{P_0^2+2P_0F_{\text{тр}}\left(1+\frac{2}{\pi}\left(1+\frac{t_0c}{2L}\right)\right)-F_{\text{тр}}^2\frac{t_0c}{2L}\right\},$$

and inversely proportional to the friction force.

The analogue formula for a rectangular pulse is:

$$Q(t) = P_0H_0(t)H_0(t_0-t). \quad (17)$$

It follows from (11) that the dependence of a velocity on time and the loading coordinate (17) is:

$$\dot{U}(z,t) = \sum_{n=0}^{n_1}\frac{1}{c\rho S_t}\left\{P_0-F_{\text{тр}}\frac{ct}{2L}\right\}H_0(t-a_{2n})H_0(a_{1n}-t)+$$
$$+\sum_{n=0}^{n_2}\frac{1}{c\rho S_t}\left\{P_0-F_{\text{тр}}\frac{ct}{2L}\right\}H_0(t-b_{2n})H_0(b_{1n}-t), \quad (18)$$

where $n_1 = \left[\frac{P_0}{F_{\text{тр}}}-\frac{z}{2L}\right]$, $n_2 = \left[\frac{P_0}{F_{\text{тр}}}+\frac{z}{2L}-1\right]$, $a_{1n} = \frac{2Ln+z}{c}+t_0$, $a_{2n} = \frac{2Ln+z}{c}$, $b_{1n} = \frac{2L(n+1)-z}{c}+t_0$, $b_{2n} = \frac{2L(n+1)-z}{c}$.

Integration of (18) in time and assumption that $z=0$ and $t > b_{1n_*}$ give the following value for the cross-section displacement $z=0$ after all the disturbances decayed in the pipe:

$$U_2^* = \frac{t_0}{c\rho S_t}\sum_{n=0}^{n_*}\left(P_0-\frac{F_{\text{тр}}c}{4L}\left(t_0+\frac{4Ln}{c}\right)\right)+\frac{t_0}{c\rho S_t}\sum_{n=0}^{n_*-1}\left(P_0-\frac{F_{\text{тр}}c}{4L}\left(t_0+\frac{4L(n+1)}{c}\right)\right), \quad n_* = \left[\frac{P_0}{F_{\text{тр}}}\right].$$

In this case all the sums are computed explicitly with the result:

$$U_2^* = \frac{t_0}{2c\rho S_t}\left\{\left(P_0 - F_{\text{тр}}\frac{t_0 c}{4L}\right)\left(1+2\left[\frac{P_0}{F_{\text{тр}}}\right]\right) - F_{\text{тр}}\left(\left[\frac{P_0}{F_{\text{тр}}}\right]+1\right)\left[\frac{P_0}{F_{\text{тр}}}\right]\right\}. \quad (19)$$

Given that $P_0 >> F_f$, (19) can be rewritten as:

$$U_2^* \approx \frac{t_0}{2c\rho S_t F_{\text{тр}}}\left\{\left(P_0 - F_{\text{тр}}\frac{t_0 c}{4L}\right)^2 - F_{\text{тр}}^2\frac{t_0 c}{4L}\right\}. \quad (20)$$

This formula is similar to (16).

## 5. COMPARISON OF NUMERICAL CALCULATION DATA WITH ANALYTICAL ESTIMATES

Figures 1 and 2 demonstrate the curves for the displacement velocity in a pipe completely embedded in the soil from the longitudinal coordinate under a semi-sinusoidal pulse without reflections $L_1 = L = 100$ m in Fig. 1 and with reflections $L_1 = L = 30$ m in Fig .2 at different time moments. The vertical dash-and-dot line corresponds to the pipe cross-section $z = z_*$. The parameters of the problem are: $P_0 = 989.6$ kN; $t_0 = 0.25$ ms; $E = 2.03\cdot 10^5$ MPa; $h = 0.01$ m; $R = 0.1625$ m; $\rho = 7805$ kg/m$^3$; $h_z = 0.1$ m; $\tau_0 = 0.02$ MPa. From the comparison of the graphs it is obvious that the numerical and analytical solutions coincide with high accuracy. It is worth noting that $t_0 / t_* = 0.013 << 1$ for the given parameters. Estimates (8) and (9) give the following values: the decaying time $t_* = 19.6$ ms for a longitudinal pulse and $z < z_* = 96.9$ m in the rod section occupied by disturbances, the cited data are confirmed by numerical calculation results in Figs. 1 and 2. The analysis of the numerical and analytical solutions in Figs. 1 and 2 reveals that the friction effect leads to the situation when a lower portion of a pulse is progressively “faded” until the pulse decays completely.

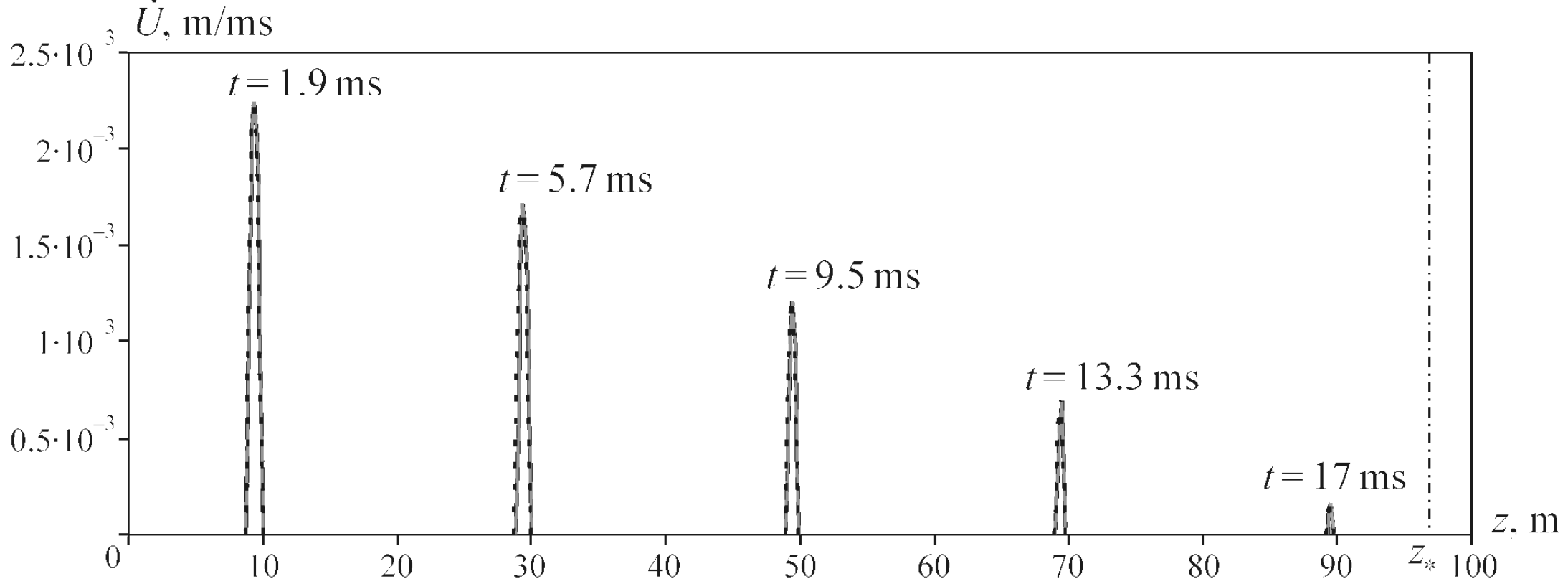

**Fig. 1.** Curves for displacement velocities at different time moments. Solid curves-numerical results of experiments, dashed curves-analytical solution (5).

In Fig. 3 the displacement curves are plotted for a pipe of 100 m in length at time moment $t = 19.6$ ms and different $\tau_0$ values. The other parameters remain the same as in Figs. 1 and 2. The analytical and numerical solutions coincide with high accuracy (Fig. 3).

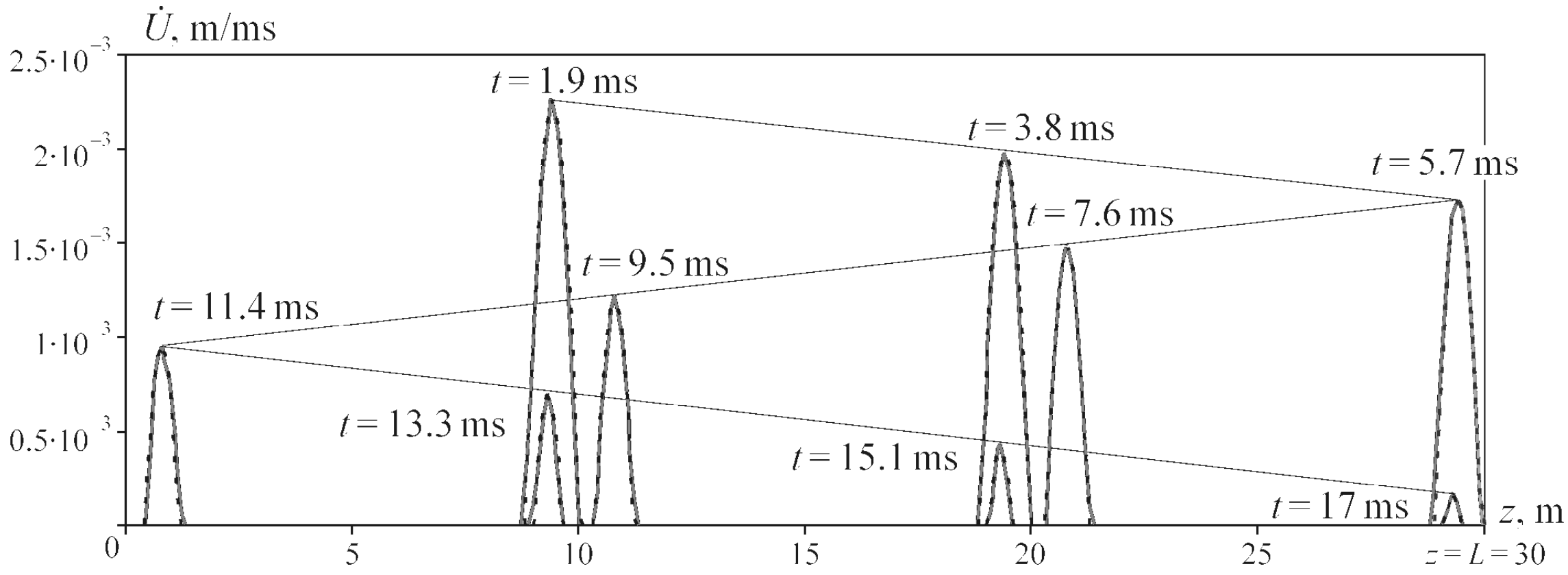

**Fig. 2.** Curves for displacement velocities at different time moments with account for reflection from rod ends at length $L = 30$ m. Solid curves-numerical results of experiments, dashed curves-analytical solution (5) and (13).

In Figs. 4 and 5 oscillograms present the numerical and analytical solutions to displacements and displacement velocities of the pipe section $z = 0$ under a semi-sinusoidal effect. The parameters of the problem are: $P_0 = 88$ kN, $t_0 = 0.22$ ms, $L = L_1 = 4$ m, $E = 2.1 \cdot 10^5$ MPa, $h = 0.003$ m, $R = 0.045$ m, $\rho = 7530$ kg/m$^3$, $h_z = 0.1$ m, $\tau_0 = 0.003$ MPa. The vertical dash-and-dot lines indicate the reflected waves at time moments $t = 2L/s$ and $t = 4L/s$. The solid line is a numerical solution and a dashed line is an analytical solution where only reflection is taken into account.

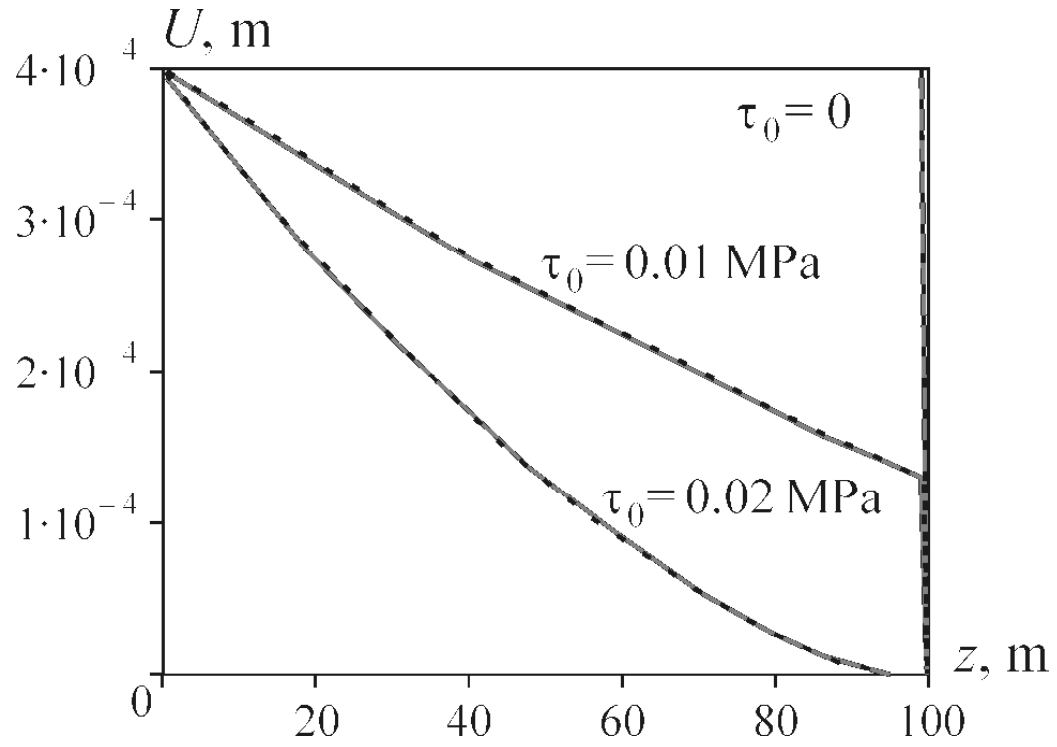

**Fig. 3.** Curves for displacements in a pipe of 100 m length at time moment $t = 19.6$ ms and different $\tau_0$ values: solid curves-numerical results of experiments, dashed curves-analytical solution (10).

The comparison of the numerical calculation results at different values of parameters in formulae (5), (10), (13), and (14) shows that the numerical and approximated analytical solutions coincide with high accuracy under terms that $t_0 << t_*$.

Figure 6 shows the numerical calculation results for displacements of a pipe section $z = 0$ versus time at multiple reflections of disturbances from the pipe ends. The parameters of the problem are identical to those in Figs. 4 and 5 at variable $\tau_0$ values. Thin lines indicate a solution under a semi-sinusoidal effect (4), Thick lines denote a solution under a rectangular pulse (17) of $t_0 / 2$ duration. The dash-and-dot lines correspond to the approximated analytical estimates (16) and (19).

Thus, a pipe slip is rather closely and qualitatively correctly evaluated by rough estimates (16) and (19) for different-shaped pulses, but with the same impact energy, amounting to $P_0^2 t_0 / 2$ in the given case, and the same friction force.

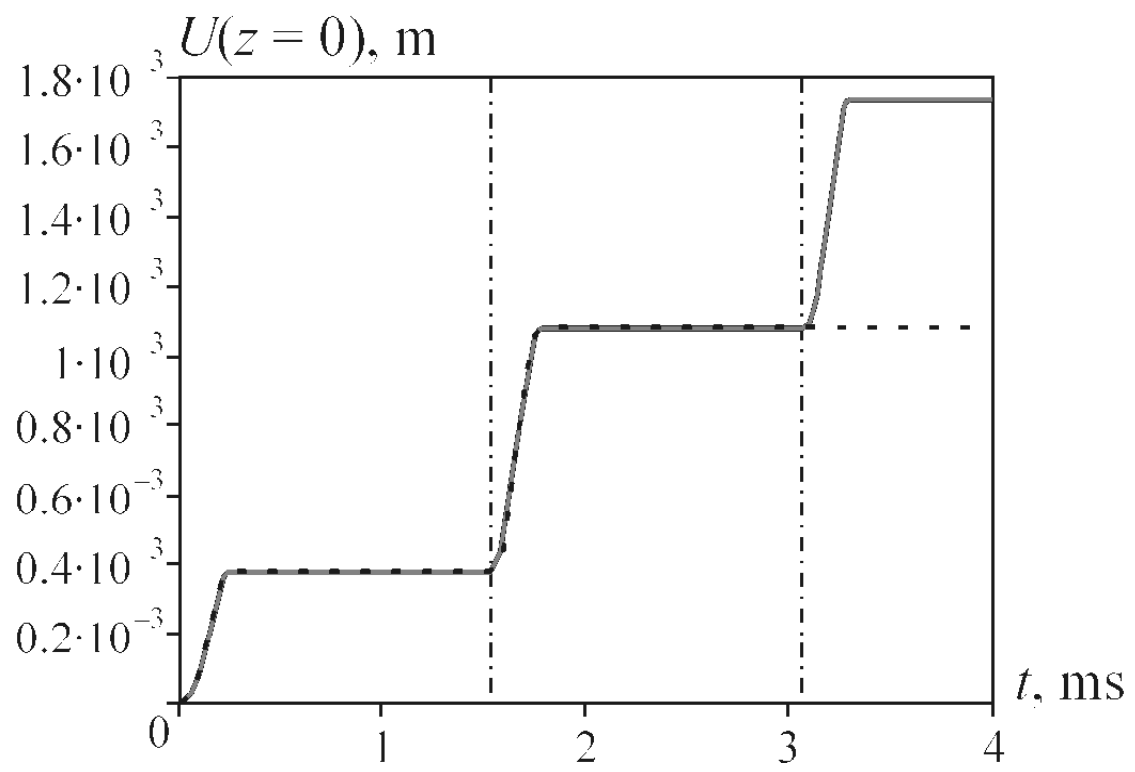

Fig. 4. Oscillograms for displacements in cross-section $z = 0$.

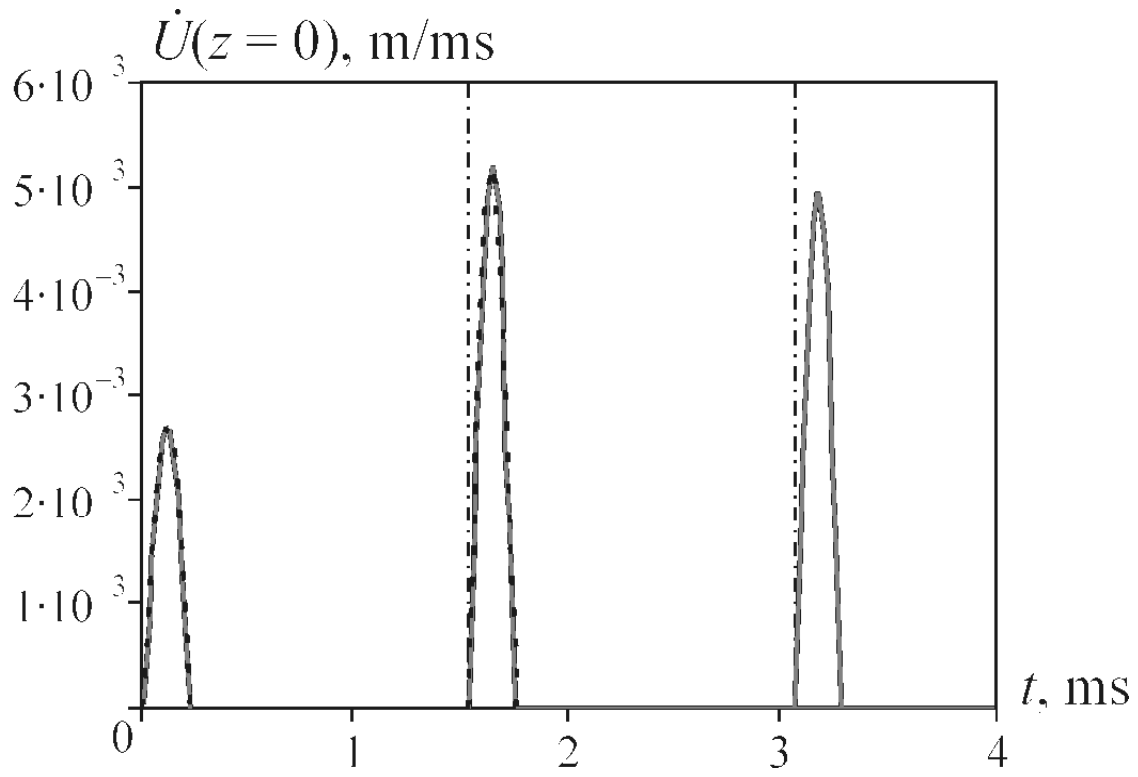

Fig. 5. Oscillograms for displacement velocities in cross-section $z = 0$.

Figure 7 presents the curves of displacements at time moment $t = 100$ ms under friction $F_f$ and semi-sinusoidal pulse. The curves are calculated by the finite difference method for a pipe partially embedded in the soil at the following parameters: $L = 7.5$ m, $L_1 = 4$ m, $P_0 = 88$ kN, $t_0 = 0.22$ ms, $h = 0.003$ m, $R = 0.045$ m, $\rho = 7530$ kg/m$^3$, $E = 2.1 \cdot 10^5$ MPa, $h_z = 0.1$ m. The solid curve is a finite difference solution, the dashed curve is an approximate approximation of a finite difference solution: $F(F_f) = 21.961 F_f^{-0.9413}$, being very close to the reverse function $21.961 / F_f$. Given that $F_f \le 2P_0$, the pipe is slipping relative to a medium. The qualitative behavior of a dependence for a pipe displacement per a single impact under the friction effect coincides with the experimental results reported in [23] for the end and cone connection of an adapter and a pipe.

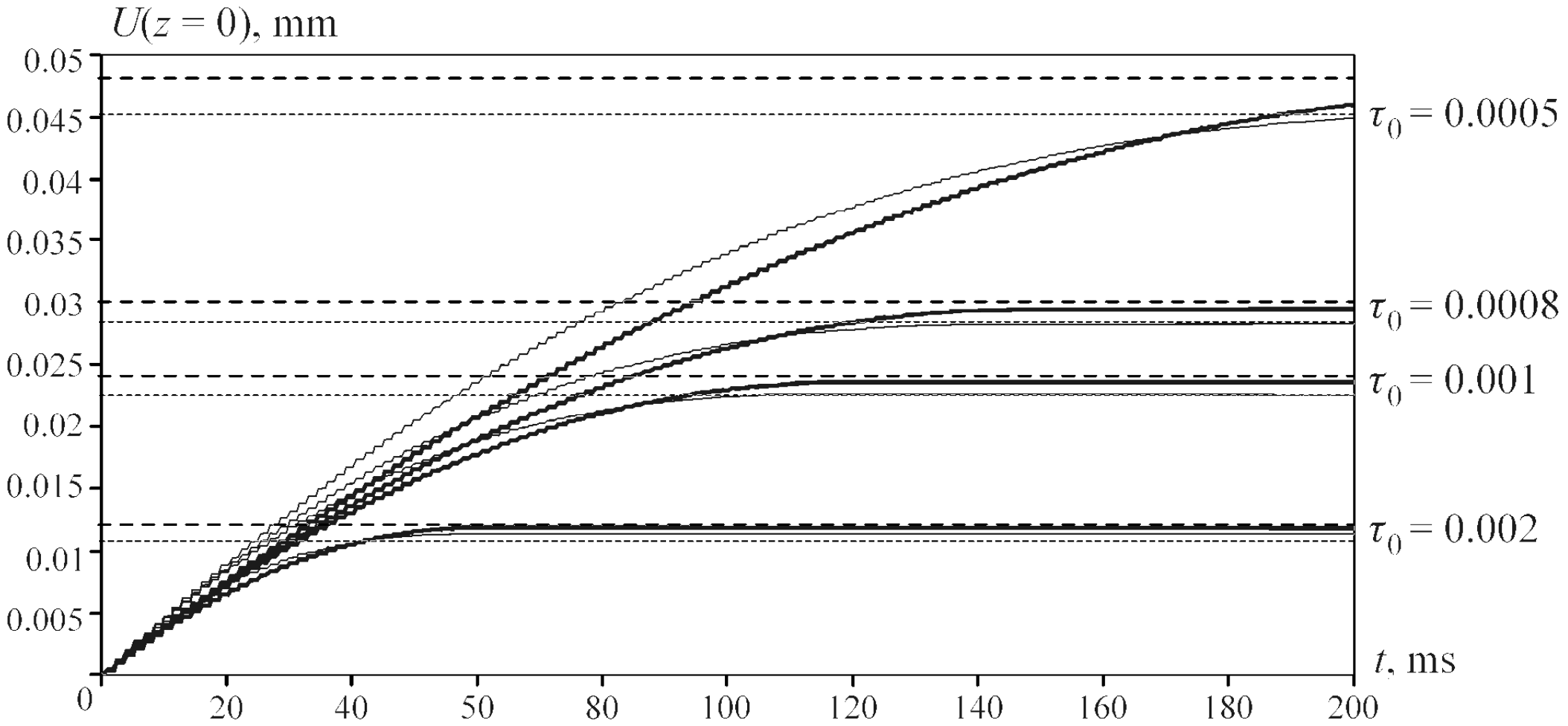

**Fig. 6.** Oscillograms for displacements of a pipe section $z = 0$: thin lines-solution for effect (4) and thick lines – for a pulse (17).

In Figs. 8 and 9 the numerical results are calculated for displacements of a pipe section $z = 0$, at different pipe length values: $L = 3.5\text{ m} + L_1$, the other parameters being identical to those in Fig. 7. Figure 8 demonstrates ocsillograms of displacements. In Fig. 9 the solid line is the dependence of displacements at time moment $t = 100$ ms on the length $L_1$ of a pipe section, embedded in the soil. The dashed line corresponds to the approximated approximation of the numerical experimental results for function $F(L_1) = 4.28 L_1^{-0.9448}$. This function is very close to the reverse function $4.28 / L_1$. Thus, for the partially-embedded-in-soil pipe the qualitative dependence on $L_1$ coincides with that obtained from formula (16) for a pipe completely embedded in the soil.

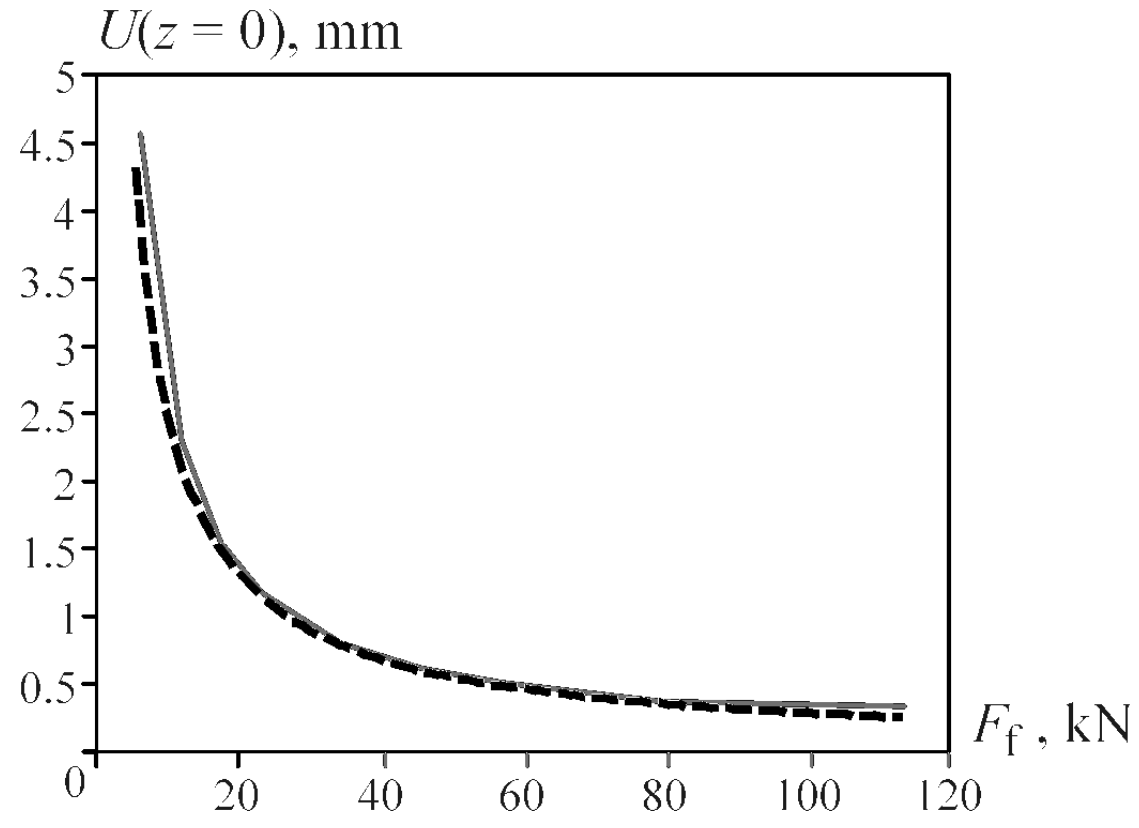

**Fig. 7.** Relationship between displacements and friction force $F_f$ at time moment $t = 100$ ms at a semi-sinusoidal pulse action.

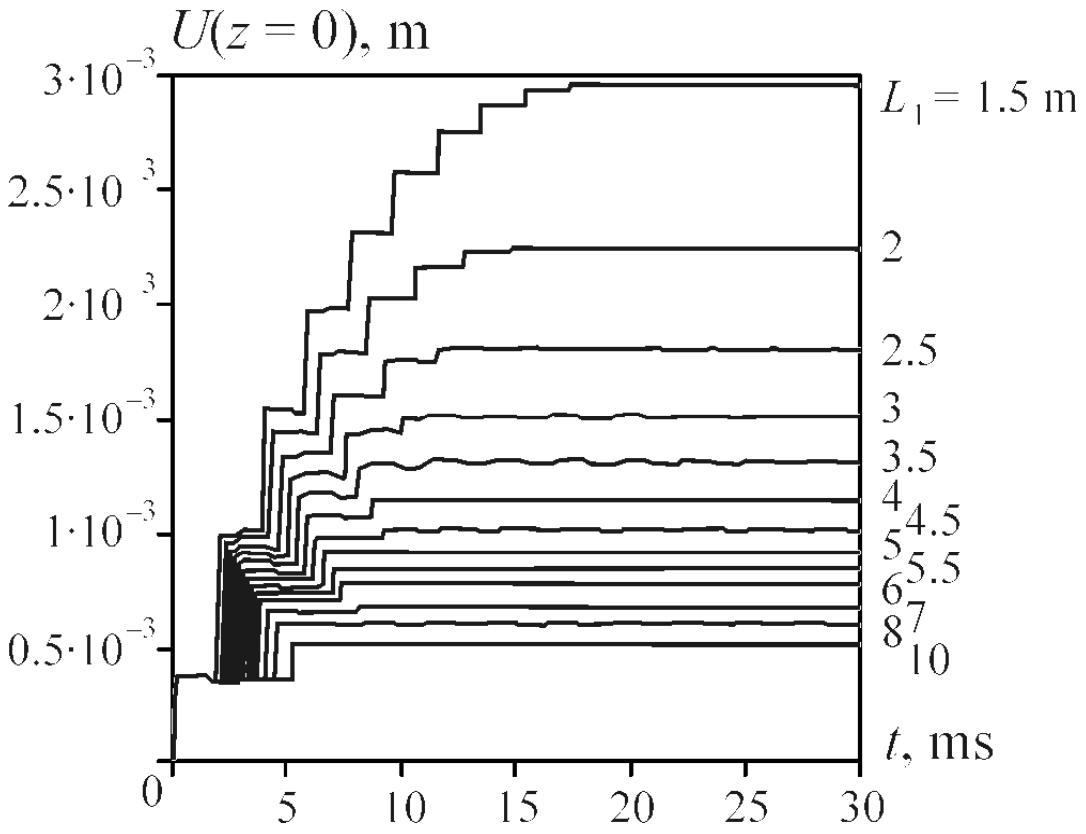

**Fig. 8.** Oscillograms for displacements at different tube lengths, embedded in the soil.

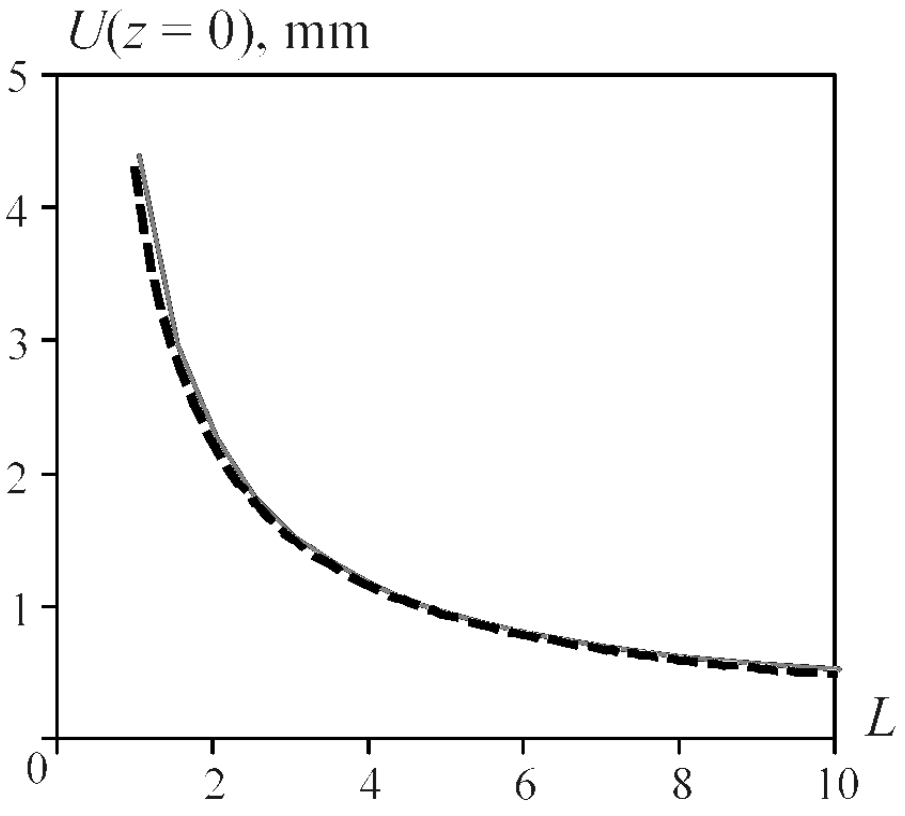

**Fig. 9.** Dependence of displacements on a length of a pipe embedded in the soil at $t = 100$ms

Dependences of the maximums of velocity on a distance obtained numerically at the parameters: $L = L_1 = 100$ m, $E = 2.03 \cdot 10^5$ MPa, $h = 0.01$ m, $R = 0.1625$ m, $\rho = 7805$ kg/m$^3$, $h_z = 0.1$ m, $\tau_0 = 0.02$ MPa, which are identical to those in [24] for the other model, are plotted in Fig. 10: curve *1* corresponds to $\sigma_0 = 50$ MPa, $t_0 = 1$ ms, curve *2* — $\sigma_0 = 70.7$ MPa, $t_0 = 0.5$ ms, curve *3* — $\sigma_0 = 100$ MPa, $t_0 = 0.25$ ms at $\sigma_0 = P_0 / S_t$. Curves *1–3* and [24] coincide qualitatively and quantitatively and are described with high accuracy by equations (5).

Figure 11 presents the distribution of displacements along the pipe length at time moment 19.6 ms for different amplitude and duration of semi-sinusoidal pulses, selected in such a way that their impact energy coincides and amounts to $P_0^2 t_0 / 2$. The solid line is an analytical solution (10), the dash-and-dot line is numerical calculations, parameters identical to those in Fig. 10. It is evident that the analytical solution describes the numerical solution with high precision except for a close-to-impact area. These curves qualitatively and quantitatively coincide with calculations obtained for the other model in [24].

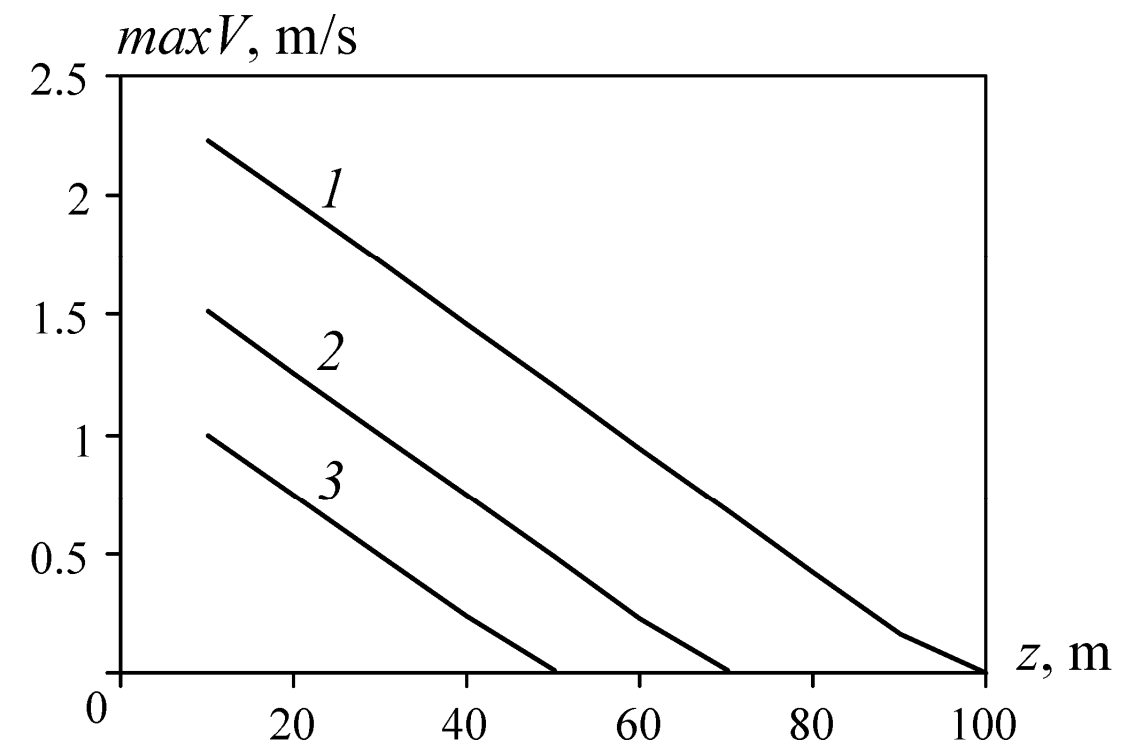

Fig. 10. Velocity maximums vs. distance (described in the text)

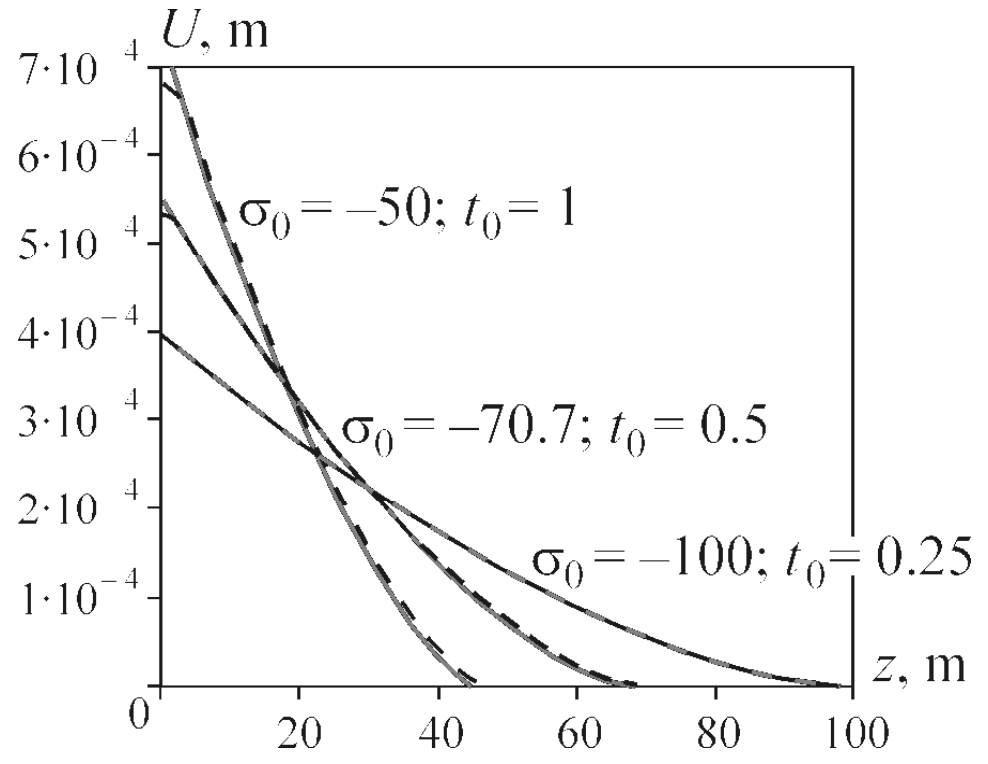

**Fig. 11**. Displacements at time moment 19.6 ms.

In Fig. 12 the results are calculated at time moment 19.6 ms under the sinusoidal load, applied to a rod end:

$$ES_t \frac{\partial u}{\partial x} = -P_0 \sin(\omega_* t) H_0(t) .$$

The parameters of the problem are: $\omega_* = 3.14$ ms$^{-1}$, $\sigma_0 = 100$ MPa, $t_0 = 1$ ms, other parameters are identical to those in Fig. 10.

In Fig. 12a the displacement curves are presented. The solid line corresponds to the finite-difference solution to equations (1)–(3), the dashed-and-dot lines is the approximated analytical solution, obtained by the method of the mechanical analogue [7]:

$$U(z,t) = \left( A_0 - \frac{2\tau_0 P_t z}{\pi c \omega_* \rho S_t} \right) \sin\left( \omega_* t - \frac{\omega_* z}{c} \right), \quad A_0 = \frac{2c}{\pi \omega_*} \sqrt{\left( \frac{\pi P_0}{2ES_t} \right)^2 - \left( \frac{\tau_0 P_t}{c \omega_* \rho S_t} \right)^2} . \quad (21)$$

The dashed line is an approximated analytical solution, obtained by the method of complex amplitudes [7]:

$$U(z,t) = \left( A_0 - \frac{\tau_0 P_t z}{2c \omega_* \rho S_t} \right) \sin\left( \omega_* t - \frac{\omega_* z}{c} \right). \quad (22)$$

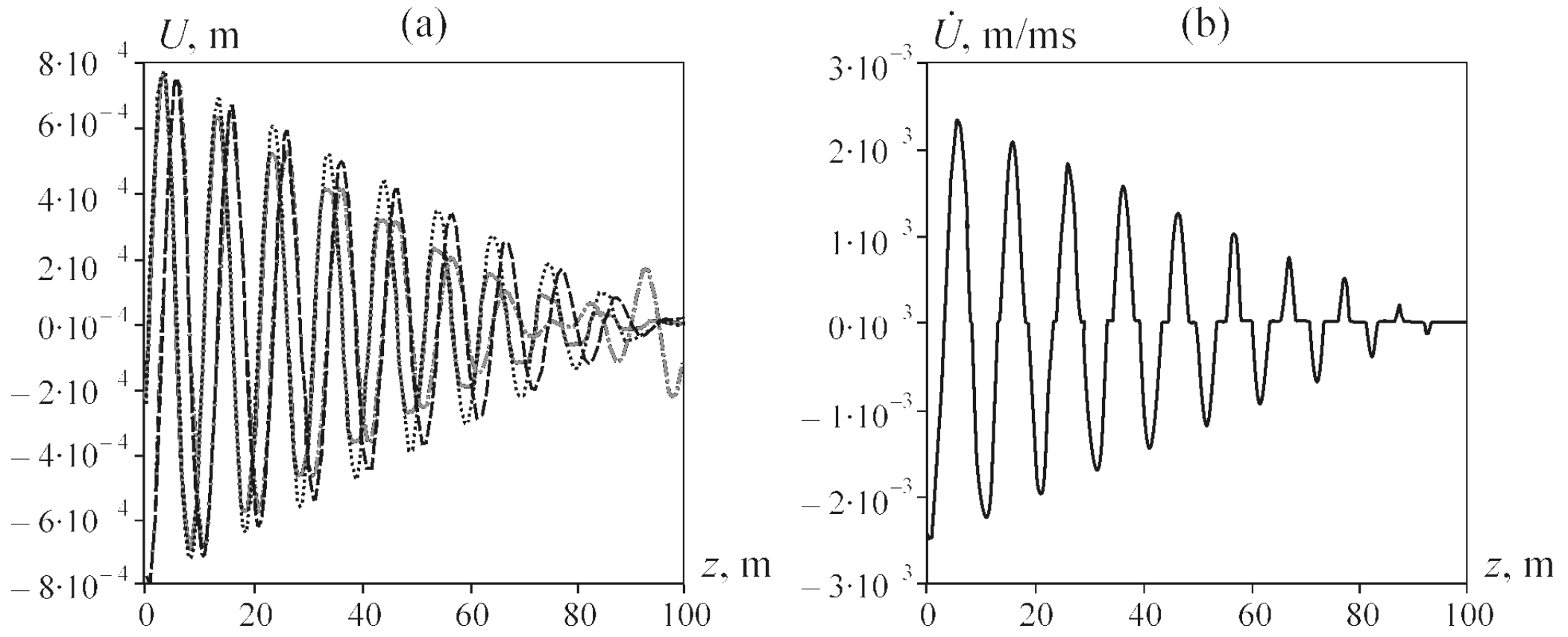

Fig.12. Displacement curves: (a) solid curve— the finite-difference solution; dashed curve—solution (22); dash-and-dot curve—solution (21); dotted curve—solution (23); (b) displacement velocity curves.

In Fig. 12bthe graphs of the displacement velocity are presented. The comparison of the cited results makes it possible to conclude that the solution obtained by the method of the mechanical analogue (21) gives a high error in an area of the running disturbance front, but with increase in the distance from the front the method gives the correct amplitude of oscillations. The solution obtained by the method of complex amplitudes (22) gives less error nearby the front, but beyond the front the oscillation amplitude is appreciably higher as compared to that, calculated in the finite-difference solution. The oscillation phase does not coincide with the numerical solution of equations (1)–(3) in both analytical solutions [7]. First, it is explicit that as the stresses at the rod end are preset by the sinusoidal function, the displacements after integration should be evaluated through cosines. Second, the second summand in the radical expression for $A_0$ amounts to less than 0.05% of the first summand at the cited parameters, so it can be neglected. Third, as it was mentioned above, the disturbances in the rod with dry friction at its side surface are absent in domain $z > z_*$. Taking into account these three remarks, we rewrite formula (22) as:

$$U(z,t) = \frac{P_0 + \tau_0 P_t / 2}{\omega_* \rho c S_t} \cos\left( \omega_* t - \frac{\omega_* z}{c} \right) H_0(z_* - z) . \tag{23}$$

In Fig. 12*a* the dashed curve corresponds to expression (23). As it follows from three solutions (21)–(23), solution (23) is the closest to the finite-difference solution. It is seen in the velocity curves how the reduction in every half-oscillation proceeds in the case of harmonic disturbances. Given that a sinusoidal effect is assumed as a sum of semi-sinusoids and consider for a velocity sign for every of these semi-sinusoids it is possible to obtain a more precise solution expressed as a sum of a series from formulae (5) and (10) as compared to the solution (23). We do not present this solution because of its bulkiness. The finite-difference solution was also compared to the well known analytical one [13], the obtained results were merging in the plots.

## CONCLUSIONS

The finite-difference algorithm for calculation of the wave processes in a pipe driven into the soil with account for the external dry friction is developed. The approximated analytical solutions to the 1D-in-longitudinal-direction problem for a pipe subjected to the applied external dry friction and a random-shaped longitudinal pulse. Solutions are obtained for both semi-infinite pipe and a finite-length pipe with account for the pulse reflections from its end. The comparison of the finite-difference and analytical solutions revealed their good compliance.

The finite-difference calculations of the disturbances in a pipe under the pulse effect with account for the contact dry friction at its side surface were made for a wide range of variations in the parameters of the problem. It is established that in the case with the finite-length pipe the pipe displacements within a large time interval, as compared to a period of to-and-fro wave run along the pipe at the velocity of longitudinal waves, is proportional to "an impact energy with account for friction" and a pulse duration and reverse proportional to the friction force.

*****